\documentclass[12pt]{article}
\usepackage[frenchb]{babel}     
\newtheorem{lem}{Lemme}[section]
\newtheorem{prop}{Proposition}[section]
\newtheorem{thm}{Th\'eor\`eme}[section]

\newtheorem{rem}{Remarque}[section]
\newtheorem{defi}{D\'efinition}[section]

\let\ch=\widehat

\let\lfd=\longrightarrow
\let\lfc=\longmapsto
\def\f#1{{\bf F}_{#1}}
\def\tr{\mathop{\rm Tr}\nolimits}
\let\dps=\displaystyle
\def\biq#1{\quad \hbox{ #1 }\quad }
\def\pt#1{\left(#1\right)}

\let\ss=\smallskip

\let\bb=\bigbreak
\def\ie{c'est-\`a-dire }
\def\mod#1{{\ ({\rm mod.}\ #1)}} 
\def\ssi{\hbox{ si et seulement si }}

\begin{document}

\title{Non lin\'earit\'e des fonctions bool\'eennes\\
donn\'ees par des traces de polyn\^omes\\ de degr\'e binaire 3}
\author{Eric F\'erard\thanks{Universit\'e de Polyn\'esie fran\c caise, Tahiti; {e-mail} {\tt ferard@upf.pf}}, 
Fran\c cois Rodier\thanks{Institut de Math\'ematiques de Luminy --
C.N.R.S. 
163 avenue de Luminy,
Case 907, Marseille Cedex 9, France; {e-mail} {\tt rodier@iml.univ-mrs.fr}; {Tel} 04 91 26 95 89;{Fax} 04 91 26 96 55
}}
\date{}
\maketitle

\begin{abstract}
Nous \'etudions la non lin\'earit\'e des fonctions d\'efinies sur $\f{2^m}$ o\`u $m$ est un entier impair, associ\'ees aux polyn\^omes de degr\'e 7 ou \`a des polyn\^omes plus g\'en\'eraux.
\end{abstract}



\section{Introduction}

La non-lin\'earit\'e d'une fonction bool\'eenne 
$f:\f2^m\lfd\f2$ 
est la distance de $f$ \`a l'ensemble des fonctions affines \`a $m$ variables (voir les \S \ \ref{defnl}). 
C'est un concept important.

Il intervient en cryptographie (cf. \cite{brv, ca2, ca3, cf}) pour construire des cryptosyst\`emes performants (chiffrements sym\'etriques), et dans la th\'eorie de codage avec le vieux probl\`eme du rayon de recouvrement des codes de  Reed-Muller d'ordre 1.

La non-lin\'earit\'e est inf\'erieure \`a $2^{m-1}-2^{m/2-1}$. Cette limite est atteinte par les fonctions courbes (cf. le livre de MacWillams et de Sloane \cite{ms}) qui existent seulement si le nombre de variables $m$ des fonctions bool\'eennes est pair. Pour des raisons de s\'ecurit\'e en cryptographie, et aussi parce que les fonctions bool\'eennes doivent avoir d'autres propri\'et\'es telles que l'\'equilibre ou le degr\'e alg\'ebrique \'elev\'e, il est important d'avoir la possibilit\'e de choix parmi beaucoup de fonctions bool\'eennes, non seulement des fonctions courbes, mais \'egalement des fonctions presque courbes dans le sens que leur non-lin\'earit\'e est voisine de la non-lin\'earit\'e des fonctions courbes.

Pour $m$ impair, il serait particuli\`erement int\'eressant de trouver des fonctions avec une non-lin\'earit\'e plus grande que celle de fonctions bool\'eennes quadratiques (appel\'ees {\sl presque optimales} dans \cite{cccf}). Ceci a \'et\'e fait dans le travail de Patterson et de Wiedemann \cite{pw} et \'egalement de Langevin et Zanotti \cite{lz} et plus r\'ecemment par Kavut, Maitra et Y{\"u}cel \cite{ma}.

\bb
Soit $q=2^m$ et $k=\f{2^m}$ assimil\'e comme espace vectoriel sur $\f2$ \`a $\f2^m$. Si $G$ est un polyn\^ome sur $k$, cela nous permet de construire une fonction bool\'eenne $\tr G(x)$, o\`u $\tr$ est la trace de $\f{2^m}$ sur $\f2$, ou plut\^ot la fonction $\chi(G(x))$, avec des valeurs dans $\pm1$, o\`u nous d\'enotons par $\chi_0$ le caract\`ere non trivial unique de $\f2$ dans les nombres complexes diff\'erents de z\'ero : $$\chi_0(0)=1\biq{,}\chi_0(1)=-1$$ 
et nous notons $\chi=\chi_0\circ\tr$.

Pour $m$ pair, on a cherch\'e \`a trouver des fonctions courbes de cette forme. Pour mentionner seulement le cas des mon\^omes, on peut consid\'erer les cas connus (de Gold, de Dillon, des exposants de Niho) dans l'article de Leander \cite{le}. Ce sont des fonctions $f:x\lfd \chi({ax^r})$ o\`u $r=3$ ou 5 (ou plus g\'en\'eralement $r=2^i+1$, o\`u $i$ est un nombre entier) et $a\in k$ n'est pas de la forme $x^r$.

On aurait pu esp\'erer que pour $r=7$, ou parmi les fonctions 
$$f:x\lfd\chi\pt{G(x)}$$
  quand $G$ est un polyn\^ome du degr\'e 7, il y a quelques fonctions qui sont presque courbes au sens pr\'ec\'edent.
Cela s'av\`ere ne pas \^etre le cas, mais nous prouverons que pour $m$ impair de telles fonctions ont les propri\'et\'es de non-lin\'earit\'e plut\^ot bonnes (cf. section \ref{deg3}).
 Nous employons pour cela des r\'esultats r\'ecents de Maisner et de Nart au sujet des fonctions de z\^eta des courbes supersinguli\`eres de genre 2 que nous avons regroup\'es dans les sections \ref{hell}, \ref{demeval}, \ref{dembornes}.

\section{Pr\'eliminaires}

\subsection{Fonctions bool\'eennes}

Soit $m$ un entier positif et $q=2^m$.

\begin{defi}
Une fonction bool\'eenne \`a $m$ variables est une application de l'espace
$V_m=(\f2)^m$ dans $\f2$.
\end{defi}

Une fonction bool\'eenne est \emph{lin\'eaire} si c'est une forme 
lin\'eaire sur
l'espace vectoriel $(\f2)^m$. Elle est dite \emph{affine} si elle est
\'egale \`a une fonction lin\'eaire \`a une constante pr\`es.

\subsection{Non-lin\'earit\'e}
\label{defnl}

\begin{defi}
On appelle  non-lin\'earit\'e d'une fonction bool\'eenne $f$ \`a
$m$ variables et on la note $nl(f)$ la distance qui la s\'epare de
l'ensemble des fonctions affines \`a $m$ variables :
$$nl(f) = \min_{h \hbox{\,\scriptsize affine }} d(f, h)$$
o\`u $d$ est la distance de Hamming.
\end{defi}

On peut prouver que la
non-lin\'earit\'e est \'egale \`a
$$\dps nl(f)  = 2^{m-1} - {1\over 2}\|\ch f\|_\infty
$$
o\`u
$$  \|\ch f\|_\infty = 
\sup_{v\in V_m} \Bigl| \sum_{x\in V_m}\chi_0{(f(x)+v\cdot x)}\Bigr|$$
et $v\cdot x$ denote le produit scalaire usuel   de $V_m$.
C'est le maximum de la transform\'ee de Fourier de  $\chi_0\pt{f}$ (ou la transform\'ee de Walsh de $f$):
$$\ch f(v) = \sum_{x\in V_m}\chi_0{(f(x)+v\cdot x)}.$$
On appellera $  \|\ch f\|_\infty$ l'amplitude spectrale de la fonction bool\'eenne $f$.
La  formule d'inversion est donn\'ee par
$$ \chi_0({f(x)}) = \frac{1}{q} \sum_{v \in V_m} \widehat{f}(v) \chi_0\pt{v.x} $$
o\`u l'on remarque que le dual de $V_m$ est isomorphe \`a $V_m$, avec la mesure  ${1\over q}$ sur chaque  point.
L'identit\'e de Parseval peut s'\'ecrire
$$\|\ch f\|_2^2= \frac{1}{q} \sum_{v \in V_m} \widehat{f}(v)^2=q$$
et,
si $f$ est une fonction bool\'eenne sur $\f2^m$:
$$\sqrt q\le  \|\ch f\|_\infty\le q.$$

\subsection{L'indice de somme des carr\'es}

Soit $f$ une fonction bool\'eenne sur $V_m$. Zhang et
 Zheng ont introduit  l'indice de somme des carr\'es   \cite{zz}:
$$ \sigma_f={1\over q}{\sum_{x\in
V_m}\ch{f}(x)^4}=\|\ch{f}\|_4^4.$$
Nous remarquons que
\begin{equation}
\label{ineg2}
\|\ch f\|_2\le\|\ch f\|_4\le\|\ch f\|_\infty.
\end{equation}
La relation de cette fonction avec la non lin\'earit\'e a \'et\'e \'etudi\'ee par A.~Canteaut et
al. \cite{cccf}.

\section{{Les fonctions $f:x\lfd \tr\pt{ G(x)}$ o\`u $G$ est un polyn\^ome}}

\subsection{Divisibilit\'e de $\|\ch f\|_\infty$}

Soit $G(x)$ le polyn\^ome $\sum_{i=0}^s a_ix^i$  \`a coefficients dans $\f q$ et $f$ la fonction bool\'eenne $\tr\circ G$.

\begin{defi}
Le degr\'e binaire de $G$ est la valeur maximum des $\sigma(i)$ pour $0\le i\le s$, o\`u
$\sigma(i)$ est la somme des chiffres de $i$ \'ecrit en chiffre binaire.
\end{defi}

On a la proposition suivante, due \`a C. Moreno et O. Moreno \cite{mm}, g\'en\'e-ralisant le th\'eor\`eme d'Ax.

\begin{prop}
Soit $G$ un polyn\^ome \`a coefficients dans $\f q$, de degr\'e binaire $d$.
Alors $\|\ch f \|_\infty$ est divisible par $2^{\lceil{m\over d}\rceil}$.
\end{prop}

\subsection{Cas o\`u $G$ est un polyn\^ome de  degr\'e binaire 2}

Les $\|\ch f\|_\infty$ sont multiples de $2^{\lceil{m\over 2}\rceil}$. Donc, si $m$ est pair $\|\ch f\|_\infty$ est un multiple de  $q^{1/2}$, et si $m$ est impair, de $\sqrt{2q}$.
En particulier, si $m$ est impair, l'amplitude spectrale est sup\'erieure ou \'egale \`a  $\sqrt{2q}$ qui est \'egale \`a celle des fonctions bool\'eennes quadratiques de rang maximal.

\section{Les fonctions $f:x\lfd \tr\pt{ G(x)}$ o\`u $G$ est un polyn\^ome de degr\'e  binaire 3}
\label{deg3}
On va simplement \'etudier le cas o\`u $G$ est un polyn\^ome de  degr\'e binaire 2 auquel on a rajout\'e un mon\^ome non nul de degr\'e 7, \ie un polyn\^ome de la forme
$$G=a_7x^7+\sum_0^s b_ix^{2^i+1}$$ 
o\`u $a_7\ne0$ un polyn\^ome de degr\'e 7 \`a coefficients dans $k$.  
Nous voudrions \'evaluer $\|\ch f\|_4$ sur $\f{2^m}$, pour $f(x)=\tr\pt{ G(x)}$ o\`u $\tr$ d\'enote la fonction trace de $\f q$ vers $\f2$:
$$\tr(x)=\sum_{i=0}^{m-1}x^{2^i}.$$

\subsection{Evaluation de $\|\ch f\|_4^4$}

\begin{prop}
\label{eval}
La valeur de $\|\ch f\|_4^4$ sur $\f{2^m}$ quand $m$ est impair et $f(x)=\chi\pt{ G(x)}$ est telle que
$$
|\|\ch f\|_4^4
-3q^2|\le 185.2^{s-1}q^{3/2}.$$
\end{prop}

{\sl D\'emonstration} --

La d\'emonstration sera donn\'e dans la section \ref{demeval}.

\begin{rem}
Ce r\'esultat est \`a comparer avec la proposition 5.6 de \cite{ro1} o\`u on a montr\'e que la distribution de $\|\ch f\|_4^4$ pour toutes les fonction bool\'eennes est concentr\'ee autour de $3q^2$.
\end{rem}

\subsection{Bornes de $\|\ch f\|_\infty$}

La d\'emonstration de ces bornes seront donn\'ees dans la section \ref{dembornes}.

\subsubsection{Borne inf\'erieure}

\begin{prop}
\label{borneinf}
Pour les fonctions $f:x\lfd \chi\pt{ G(x)}$ sur $\f{2^m}$ o\`u $G$ est un polyn\^ome donn\'e au d\'ebut de la section \ref{deg3} et $m$ est impair, on a, pour $m\le11+2s$
$$\sqrt{2 q} \le \|\ch f\|_\infty.$$
Pour $m\ge15+2s$, on a de plus
$$\sqrt{2 q}+2^{\lceil{m\over 3}\rceil}\le \|\ch f\|_\infty.$$
\end{prop}

\begin{rem}
\label{connu}
Il est connu que pour $m$ impair et plus petit que 7, on a
$\sqrt{2 q}\le \|\ch f\|_\infty$
pour toutes les  fonctions bool\'eennes.
\cite{ho}
\end{rem}

\subsubsection{Borne sup\'erieure}

\begin{prop}
\label{bornesup}
On a
$$\|\ch f\|_\infty\le6\sqrt q.$$
\end{prop}

\section{Etude de courbes hyperelliptiques}
\label{hell}

Pour d\'emontrer les r\'esultats pr\'ec\'edents, on va \'etudier des courbes li\'ees au polyn\^ome $G$.

On obtient d'abord l'expression simple de $\|\ch f\|_4$ (cf. \cite{ro1,ro2}):
\label{som}
$$\|\ch f\|_4^4 = \sum_{x_1+x_2+x_3+x_4=0}
\chi\pt{f(x_1)+f(x_2)+f(x_3)+f(x_4)}=q^2+\sum_{\textstyle{\alpha\ne0\atop \alpha\in V_m}}X_\alpha$$
avec
$$ X_\alpha=\Big(\sum_{x\in k}\chi\pt{ G(x)+ G(x+\alpha)}\Big)^2.$$
On note maintenant $\alpha$ un \'el\'ement de 
$k^*$.
On peut v\'erifier que
$$\displaylines{
G(x+\alpha)+G(x)=G(\alpha)  + a_7 \alpha^6 x  + a_7 \alpha^5 x^2 + 
    a_7 \alpha^4 x^3  + a_7 \alpha^3 x^4 + a_7 \alpha^2 x^5 + \cr
    a_7 \alpha x^6 +\sum_0^s b_i(\alpha x^{2^i}+\alpha^{2^i} x)
}$$
Pour calculer $X_\alpha$, on peut remarquer que
la courbe d'\'equation $y^2+y=G(x+\alpha)+G(x)$ est isomorphe \`a
\begin{eqnarray*}
\lefteqn{y^2+y =
G(\alpha)+ }\\
    &&+ \Big(a_7 \alpha^6  
     + a_7^{1/4} \alpha^{3/4}  + a_7^{1/2}  \alpha^{5/2}
     +\sum_0^s (b_i\alpha)^{2^{-i}} +\sum_0^s b_i\alpha^{2^i} \Big) x +\cr 
 &&\quad  
   +( a_7 \alpha^4+a_7^{1/2} \alpha^{1/2}) x^3 + a_7 \alpha^2 x^5 
\end{eqnarray*}
qui est une \'equation de la courbe $C_1$ de genre 2 pour $\alpha\ne0$.

On a
$$X_\alpha=(\#C_1-q-1)^2.$$

\subsection{La th\'eorie de van der Geer et van der Vlugt}

Soit $C_1$ la courbe d'\'equation affine:
$$C_1 : y^2 + y = ax^5 + bx^3 + cx + d$$
avec $a\ne0$.
Soit $R$ le polyn\^ome lin\'eaire $ax^4 + bx^2 + c^2x$.
L'application 
\begin{eqnarray*}
Q : k &\lfd& \f2\\ 
x &\lfc&  \tr(xR(x))
\end{eqnarray*}
 est la forme quadratique associ\'e \`a la forme symplectique
\begin{eqnarray*}
k \times k &\lfd& \f2\\ 
(x, y) &\lfc& < x, y >= \tr(xR(y) + yR(x)).
\end{eqnarray*}
Le nombre de z\'eros de $Q$ d\'etermine le nombre de points de $C_1$:
$$\#C_1(k)=1+2\#Q^{-1}(0)$$
Le radical $W$ de la forme symplectique $<,>$ co•ncide avec l'ensemble des z\'eros dans $k$
du polyn\^ome $\f2$-lin\'eaire et s\'eparable
$$E_{a,b} = a^4x^{16} + b^4x^8 + b^2x^2 + ax.$$
On a : 
$0 \le w = dim_{\f2} W \le 4$ et $w \equiv m \mod 2$. La codimension du noyau $V$ de $Q$ dans $W$ est \'egale \`a 0 ou 1.
De plus, le polyn\^ome $E_{a,b}$ se factorise dans $k[x]$ (\cite{gv1}, Theorem 3.4):
$$E_{a,b}(x)  = xP(x)(1 + x^5P(x))$$
avec $ P (x) = a^2x^5 + b^2x + a.$

\begin{thm}
(van der Geer - van der Vlugt \cite{gv1}) 
\label{gv}

Si $V \subset W$, alors
$\#C_1(k) = 1 + q.$

Si $V = W$, alors
$\#C_1(k) = 1 + q \pm \sqrt{2^wq}.$

\end{thm} 

\subsection{Les travaux de Maisner et Nart}

Supposons que $a=b$ et que le polyn\^ome $P$ ait au moins une racine $z$.
Alors, comme $m$ est impair, il existe un unique $\ell\in k$ tel que $\ell^3=1+z^{-4}$.  

\begin{prop}
\label{tracel}
Si $\tr\ell=0$ alors le polyn\^ome $P$ a exactement trois racines dans $k$ et on a $w=3$.
Si $\tr\ell\ne0$ alors le polyn\^ome $P$ n'a qu'une  racine dans $k$, la composante restante est irr\'eductible et on a $w=1$.
\end{prop}

{D\'emonstration} --

Voir Maisner et Nart \cite{mn} propositions 2.3 et 2.6.

\subsection{R\'eduction de la courbe $y^2+y=G(x+\alpha)+G(x)$}

Soit
$\lambda=\alpha+a_7^{-1/4}\alpha^{-3/4}$.

\subsubsection{Cas o\`u $\lambda=0$}
\label{lambda0}

Alors on a
$\alpha^7=a_7^{-1}$, donc l'\'equation de la courbe devient
\begin{eqnarray*}
y^2+y&=&
G(\alpha)
    + (a_7 \alpha^6  
     + a_7^{1/4} \alpha^{3/4}  + a_7^{1/2}  \alpha^{5/2}
     +\sum_0^s b_i(\alpha^{2^{-i}} +\alpha^{2^i} )) x +\cr 
 &&\quad  
    + a_7 \alpha^2 x^5 \\
&=&d+c x +ax^5
\end{eqnarray*}
pour  $a= \alpha ^{-5}$.

    Le polyn\^ome $P$  s'\'ecrit
    $P(x)=a^2x^5+a$.
    Si $m$ est impair
  il a une unique racine
    $z=a^{-1/5}=\alpha$.
    D'apr\`es Maisner et Nart (\cite{mn}, Propositions 2.5 et 2.3)
on est dans le cas o\`u $w=1$ donc $W=\{0,z\}$.
Soit $c$ le coefficient de $x$.
On a
\begin{eqnarray*}
\tr (cz)&=&\tr \Big((1/ \alpha 
     + a_7^{1/4} \alpha^{3/4}  + a_7^{1/2}  \alpha^{5/2}
         +\sum_0^s (b_i\alpha)^{2^{-i}} +\sum_0^s b_i\alpha^{2^i})\alpha\Big)\\
    &=& \tr(1  +\sum_0^s b_i^{2^{-i}} \alpha^{2^i+1\over2^{i}}  +\sum_0^s b_i\alpha^{1+2^i} )\\
    &=&\tr 1=1
\end{eqnarray*}
On v\'erifie alors que
$$Q(z)=\tr{(az^5+cz)}=\tr{(1+cz)}=0$$
D'o\`u $V=W$ et donc
$X_\alpha=2q$ par le th\'eor\`eme \ref{gv}.

\subsubsection{Cas o\`u $\lambda\ne0$}

Cette courbe est isomorphe \`a
$$y^2+y=ax^5+ax^3+cx+d$$
avec
$$a=\lambda^5a_7 \alpha^2=\lambda^3( a_7 \alpha^4+a_7^{1/2} \alpha^{1/2})
$$
et
$\lambda=\alpha+a_7^{-1/4}\alpha^{-3/4}$.
On a
\begin{equation}
\label{calcula}
a=1 + a_7^{-1/4} \alpha^{-7/4} + a_7^{3/4}  \alpha^{21/4} + 
    a_7  \alpha^7
\end{equation}
    et
\begin{eqnarray}
\label{calculc}
c&=&  1+\Big(\sum_0^s (b_i\alpha)^{2^{-i}} +\sum_0^s b_i\alpha^{2^i}\Big)\lambda +
{a_7}^{1/2}  \alpha  ^{7/2 } + 
    a_7 ^ {3/4 }  \alpha  ^ {21/4 } + a_7  \alpha  ^7.
\end{eqnarray}

    \subsection{Valeurs de $X_\alpha$}

 \begin{prop}
 \label{Xalpha}
Supposons que $m$ soit impair. Alors
$
 X_\alpha=0 \hbox{ , } 2q  \hbox{ ou } 8q.$
\ss
Soit
$\ell=a_7 ^{-1/3}\alpha^{-7/3}$.
Alors
\begin{eqnarray}
 X_\alpha=8q&&\hbox{\ssi}\hfill\nonumber\\
&&\tr\ell=0 \biq{,} \ell=v+v^4  
\biq{,} 
\nonumber\\\hfill
&&\tr\pt{{\eta} v^3}=1 \biq{,}
\tr\pt{{\eta} (v+v^2)}=1\quad;\nonumber\\
&&\mbox{avec}\quad \eta=  1  +\sum_0^s (b_i\alpha^{1+2^{i}} )^{2^{-i}}+\sum_0^s b_i\alpha^{1+2^i}  +
\nonumber\\&&\hskip40mm+ 
{a_7}^{1/2}  \alpha ^ {7/2 } + a_7 ^ {1/4 }  \alpha ^ {7/4 }; \label{eta}\\
\ X_\alpha=2q&&\ssi
\tr\ell=1 \biq{;} \hfill\nonumber\\
 X_\alpha=0&&
\hbox{dans les  cas restant.}\hfil\nonumber
\end{eqnarray}
\end{prop}

{\sl D\'emonstration} --

Si $\lambda=0$, alors $\ell=1$ d'o\`u $\tr\ell=1$. On a bien $X_\alpha=2q$ d'apr\`es \ref{lambda0}. 

Si $\lambda\ne0$,
on \'etudie le polyn\^ome
$P=a^2x^5+a^2x+a.$
Remarquons que $z=\lambda^{-1}\alpha$ est racine de $P$.
Donc
\begin{eqnarray*}
P&=&(x+z)( a^2 x^4 + a^2 x^3 z + a^2 x^2 z^2 + a^2 x z^3 + a^2 z^4+a^2 )\cr
&=&a^2z^{-4}(x+z)(  {x^4 z^{-4}}+  x^3 z^{-3} +  x^2 z^{-2} +  x z^{-1} +  z^{-4}+ 1).
\end{eqnarray*}
La d\'ecomposition de $P$ en composante irr\'eductibles d\'epend de $e=1+z^{-4}$.
On a
$$e=1+z^{-4}=1+\lambda^4\alpha^{-4}=1+(\alpha^4+a_7 ^{-1}\alpha^{-3})\alpha^{-4}=
1+(1+a_7 ^{-1}\alpha^{-7})=a_7 ^{-1}\alpha^{-7}.$$
Comme $m$ est impair, on a $k^3=k$. Soit $\ell=e^{1/3}$. Alors, d'apr\`es la proposition \ref{tracel}, on a
\[
\left\{
\begin{array}{rcl}
 w=1 &\hbox{si} & \tr\ell=1  \\
 w=3 & \hbox{si}  & \tr\ell=0  
 \end{array}
\right.
\]
D'apr\`es le th\'eor\`eme \ref{gv}, on a

dans le premier cas, 
$X_\alpha=0\biq{ou}2q.$

dans le deuxi\`eme cas, 
$X_\alpha=0\biq{ou}8q.$
\bb
\paragraph{Premier cas, $\tr\ell=1$.}

On a $W=\{0,z\}$
et 
$$Q(z)=\tr(az^5+az^3+cz)
=\tr(az+cz+1)$$
car $\tr(az^3)=0$.
Pour que 
$X_\alpha=0$
il faut et il suffit que  $\tr (a+c)z=0$.
Des \'equations (\ref{calcula}) et (\ref{calculc}) on d\'eduit
\begin{eqnarray*}
(a+c)z
&=&
 1 + a_7 ^ {1/4 }  \alpha ^ {7/4 } +  a_7^{1/2}  \alpha ^ {7/2 } +\Big(\sum_0^s {b_i\alpha}^{2^{-i}} +\sum_0^s b_i\alpha^{2^i}\Big)\alpha   \cr
&=&
 1 + a_7 ^ {1/4 }  \alpha ^ {7/4 } +  a_7^{1/2}  \alpha ^ {7/2 } +\Big(\sum_0^s (b_i\alpha^{1+2^{i}} )^{2^{-i}}+\sum_0^s b_i\alpha^{1+2^i}\Big) 
\end{eqnarray*}
Donc
$\tr ((a+c)z)=\tr 1=1$
et
$X_\alpha=2q$.

\paragraph{Deuxi\`eme cas, $\tr\ell=0$.}

On a $W=<z,z_1,z_2>$.

Pour que 
$X_\alpha=0$
il faut et il suffit que $\tr (a+c)z_i=0$ pour l'un des $i=1,2$ ou que $\tr (a+c)z=0$.

Les nombres $z_i$ sont racines de $ {x^4 z^{-4}}+  x^3 z^{-3} +  x^2 z^{-2} +  x z^{-1} +  z^{-4}+ 1=0$.
On a $e=1+z^{-4}=\ell^3$ et $\ell=u+u^2$. D'o\`u, d'apr\`es Maisner et Nart \cite{mn} (d\'emonstration du lemme 2.4):
$$\displaylines{
{x^4 z^{-4}}+  x^3 z^{-3} +  x^2 z^{-2} +  x z^{-1} +  z^{-4}+ 1
=\hfill\cr
(x^2 z^{-2} +u x z^{-1}+(1+u)^3)(x^2 z^{-2} +(u+1) x z^{-1}+u^3)
}$$
On peut supposer $\tr u=0$ (car $\tr 1=1$, donc $u$ ou $1+u$ a une trace nulle). Soit donc $u=v+v^2$. On a par cons\'equent
$\ell=v+v^4$.
Alors le polyn\^ome
$x^2 z^{-2} +(u+1) x z^{-1}+u^3$ est r\'eductible:
ses racines sont:
$z(v(1+u)+1)=z(v(1+v+v^2)+1)=z(v^3+v+v^2+1)$ et $z(v(1+u)+u)=zv^3$.

\subsection{Calcul du nombre des $\alpha$ donn\'es par la proposition \ref{Xalpha}}

On peut \'evaluer le nombre des $\alpha$ qui donnent chaque cas de la  proposition pr\'ec\'edente.

\subsubsection{Le nombre des $\alpha$ tels que $ X_\alpha=2q$}

D'abord,  on \'evalue le nombre des $\alpha$ tels que $\tr \ell=1$ dans la proposition \ref{Xalpha}.

\begin{prop}
Le nombre $N_0$ des valeurs de $\alpha$ telles que $X_\alpha=2q$ v\'erifie
$$\Big|N_0-{q\over2}\Big|<3q^{1/2}.$$
\end{prop}

{\sl D\'emonstration} --

On a
$\tr\ell=\tr(a_7 ^{-1/3}\alpha^{-7/3})$.
Le nombre de $\alpha$ dans $k^*$ tels que $\tr(a_7 ^{-1/3}\alpha^{-7/3})=1$ est 
\'egal au nombre $N_0$ de $x$ dans $k^*$ tels que 
$\tr(a_7 ^{-1/3}x^7)=1$.
D\'efinissons
$$S_0=\sum_{x\in k} \chi{(a_7 ^{-1/3}x^7)}=N_0-(q-N_0)=2N_0-q.$$
On a
$|S_0|< 6\sqrt q$ d'o\`u
$${q-6\sqrt q\over2}\le N_0={S_0+q\over2}\le {q+6\sqrt q\over2}.$$

\subsubsection{Une courbe auxiliaire}

On a besoin d'\'evaluer le nombre des $(\alpha,v)$ v\'erifiant certaine conditions, avec $v$ tel que $v+v^4=\ell=a_7 ^{-1/3}\alpha^{-7/3}$.
Soit $x^{-3}=\alpha$
et
$a_7 ^{-1/3}=\gamma$.

\begin{prop}
On consid\`ere la courbe $C$ donn\'ee par l'\'equation
$$v+v^4=\gamma x^7$$
avec les coordonn\'ees $x$ et $v$ et le mod\`ele non singulier $\tilde C$.
Alors le morphisme $\tilde C\lfd C$ est bijectif. La courbe a un unique point \`a l'infini.
Elle est de genre 9.
Les valuations au point $(0,0)$ sont
$v_{(0,0)}(x)=1$ et $v_{(0,0)}(v)=7$.
Les valuations au point  \`a l'infini sont
$v_\infty(x)=-4$ et $v_\infty(v)=-7$.
\end{prop}

{\sl D\'emonstration} --

Voir le livre de Stichtenoth \cite{sti} p. 200.

\subsubsection{Bornes pour les sommes exponentielles}

Sur la courbe $\tilde C$, on consid\`ere une fonction rationnelle $f$, {qui n'est pas de la forme $\phi^2+\phi$, avec $\phi$ un fonction rationnelle sur $\tilde C$}. 
Soit
$$S=\sum_{z\in \tilde C_0(k)}{}'\chi\pt{ f(z)}$$
o\`u la somme est d\'efinie sur les points rationnels sur $k$ de $\tilde C$, qui ne sont pas des p\^oles de $f$.
Soit $(f)_\infty$ le diviseur des p\^oles de $f$ et $t$ le nombre de p\^oles de $f$, sans multiplicit\'e. La proposition suivante donne une  borne pour les sommes exponentielles $S$.

\begin{prop}
On a
$$|S|\le(2 g-2+t+\deg (f)_\infty)\sqrt q$$
\end{prop}

{\sl D\'emonstration} --

Voir  l'article de Bombieri, \cite{bo}.

\subsubsection{Le nombre des $(\alpha,v)$ tels que $\tr\pt{{\eta} v^3}=1$}

On \'evalue le nombre des $(\alpha,v)$ tels que 
$\tr\pt{{\eta} v^3}=1$, o\`u $\eta$ est donn\'e par  (\ref{eta}).

On a
\begin{eqnarray*}
\tr\pt{{\eta} v^3}&=&
\tr\pt{v^3  +\sum_0^s v^3(b_i\alpha^{1+2^{i}} )^{2^{-i}}+\sum_0^s v^3b_i\alpha^{1+2^i}  +  v^3\sqrt{a_7}  \alpha ^ {7/2 } + v^3a_7 ^ {1/4 }  \alpha ^ {7/4 }}\\
&=&\tr\pt{v^3  +\sum_0^s (v^{3.2^i}+v^3)b_i x^{-3-3.2^{i}}   +  (v^6 + v^{12})a_7  x ^ {-21 }}.
\end{eqnarray*}
Sur la courbe $C$, on consid\`ere la fonction
\begin{eqnarray*}
f(x)&=&v^3  +\sum_0^s (v^{3.2^i}+v^3)b_i x^{-3-3.2^{i}}   +  (v^6 + v^{12})a_7  x ^ {-21 }.
\end{eqnarray*}
Pour v\'erifier que $f$ n'est pas de la forme $\phi^2+\phi$, on consid\`ere
 $$\psi=\gamma^{1/4}{x\over v}(v^3x^{-3})^{2^{i-2}}=\gamma^{1/4}(vx^{-1})^{3.2^{i-2}-1}.$$
 Si $i\ge2$, on a
\begin{eqnarray*}
(v^{3.2^i}+v^3)x^{-3-3.2^{i}} +\psi^4+\psi
&=&\Big( {x^{-3}(\gamma x^7+v)+\gamma x^4\over v^4}\Big)(v^3x^{-3})^{2^{i}}+v^3x^{-3-3.2^{i}}+\psi\\
&=&( {x^{-3} v^{-3}})(v^3x^{-3})^{2^{i}}+v^3x^{-3-3.2^{i}}+\psi
\end{eqnarray*}
Et sa valuation \`a l'infini est donn\'ee par
\begin{eqnarray*}
v_\infty\Big((v^{3.2^i}+v^3)x^{-3-3.2^{i}} +\psi^4+\psi\Big)
&=&v_\infty\Big(( {x^{-3} v^{-3}})(v^3x^{-3})^{2^{i}}+v^3x^{-3-3.2^{i}}+\psi\Big)\\
&=&33-9.2^i
\end{eqnarray*}
si $i\ge3$.
C'est un entier  n\'egatif impair.

En faisant de m\^eme pour chaque entier $i$ dans l'expression de $f$, on trouve une fonction $\psi$ telle que la valuation au point \`a l'infini de $f+\psi^2+\psi$ soit un entier impair n\'egatif.

On peut v\'erifier que la fonction $f$ est d\'efinie sur chaque point fini de $C$ sauf peut-\^etre en les points tels que $x=0$.

On consid\`ere la somme
$$S_1=\sum_{(x,v)\in C(k)-C_\infty}\chi\pt{ f}$$
o\`u
$C_\infty=\{(0,0),(0,1),(0,\beta),(0,\beta^2),\infty\}$
et $\beta$ est une racine primitive $3^{\hbox{\tiny \`eme}}$ de l'unit\'e.
Les p\^oles de $f$ ne peuvent \^etre que parmi les points dans $C_\infty$.
La valuation de $f$ \`a l'infini est
$$v_\infty(f)\ge\inf(v_\infty(v^3),v_\infty(b_s   x ^ {-3(1+2^s) } v^{3.2^s}))\ge-9.2^s+12$$
si $s\ge2$.
La valuation de $f$ en $(0,0)$ est
$$v_{(0,0)}(f)=v_{(0,0)}(v^{3} x^{-3-3.2^{s}} )=21-3(1+2^s)=18-3.2^s.$$
La valuation de $v^{3.2^i}+v^3$ en $(0,1)$ est
\begin{eqnarray*}
v_{(0,1)}(v^{3.2^i}+v^3)&=&v_{(0,1)}\Big((v^3+1)\prod _{\delta\in\f{2^i}-\{1\}}(v^3-\delta)\Big)
=v_{(0,1)}\Big({x^7\over v}\Big)=7.
\end{eqnarray*}
La valuation de $(v^{3.2^i}+v^3)x^{-3-3.2^{i}}$ en $(0,1)$ est donc
$$ v_{(0,1)}(v^{3.2^i}+v^3)x^{-3-3.2^{i}}=7-3-3.2^{i}=4-3.2^{i}.$$
La valuation de $v^6+v^{12}$ en $(0,1)$ est
\begin{eqnarray*}
v_{(0,1)}(v^6+v^{12})=2v_{(0,1)}(1+v^3)
=2v_{(0,1)}(x^7)=14.
\end{eqnarray*}
La valuation de $(v^6+v^{12})x^{-21}$ en $(0,1)$ est
$$ v_{(0,1)}((v^6+v^{12})x^{-21})=14-21=-7.$$
La valuation de $f$ en $(0,1)$ est finalement
$$v_{(0,1)}(f)
=\inf(4-3.2^{s},-7)=4-3.2^{s}$$
si $4-3.2^s<-7$ \ie si $s\ge2$.

Le m\^eme calcul vaut pour la valuation de $f$ en $(0,\beta)$.
On a donc, pour $s\ge2$:
$$\deg(f)_\infty=-3(4-3.2^s)-(18-3.2^s)-12+9.2^s=-42 + 21 .2^s$$
et
$$|S_1|\le (18-2+5-42 + 21 .2^s) q^{1/2}=(21.2^s-21) q^{1/2}.$$

Consid\'erons sur la courbe $C$
le nombre $N_1$ des couples $(\alpha,v)$ tels que $\tr\pt{{\eta} v^3}=1 $.
Alors
{\begin{eqnarray*}
S_1=\sum_{(x,v)\in C-C_\infty} \chi\pt{ f}
=\sum_{\tr f=0} 1-N_1
=\#C-2N_1-5
\end{eqnarray*}}
 o\`u $\#C$ est le nombre des points de la courbe $C$. Donc 
$$\Big|N_1-{\#C\over2}\Big|={|S_1+5|\over2}\le{21.2^s-21\over2} q^{1/2}+5/2.$$

\subsubsection{Le nombre des $(\alpha,v)$ tels que $\tr\pt{{\eta} (v^2+v)}=1$}
\label{2cas}

Ensuite, nous \'evaluons le nombre des $(\alpha,v)$ tels que 
$\tr\pt{{\eta} (v^2+v)}=1$, o\`u $\eta$ est donn\'e par  (\ref{eta}).
\begin{eqnarray*}
\tr\pt{{\eta} (v^2+v)}&=&
\tr\Big( (a_7 ^ {1/4 }  \alpha ^ {7/4 } +  a_7^{1/2}  \alpha ^ {7/2 } +(\sum_0^s (b_i\alpha^{1+2^{i}} )^{2^{-i}}+\sum_0^s b_i\alpha^{1+2^i})
)(v^2+v)\Big) \\
&=&
\tr\Big( a_7 \gamma^2 x^ {-7 } 
+ \sum_0^s (b_i x^{-3(1+2^{i})} )(v^{2^{i+1}}+v^{2^{i}}+v^2+v)\Big) .
   \end{eqnarray*}
   On d\'efinit la
fonction $g(x)= a_7 \gamma^2 x^ {-7 } 
+ \sum_0^s (b_i x^{-3(1+2^{i})} )(v^{2^{i+1}}+v^{2^{i}}+v^2+v) $.
Elle n'est pas de la forme $\phi^2+\phi$ parce que
avec $\psi=\gamma^{1/2}x^{2-3.2^{s-1}} $, on a
\begin{eqnarray*}
v_{0,0}(x^{-3(1+2^{s})} v+\psi^2+\psi)
&=&v_{0,0}(x^{-3(1+2^{s})} v+\gamma(x^{4-3.2^{s}})+\gamma^{1/2}x^{2-3.2^{s-1}})\\
&=&v_{0,0}(x^{-3.2^s} (x^{-3}(v^4+\gamma x^7)+\gamma x^{4})+\gamma^{1/2}x^{2-3.2^{s-1}})\\
&\ge&\inf(-3.2^s+25,2-3.2^{s-1})
\end{eqnarray*}
d'o\`u
\begin{eqnarray*}
v_{0,0}(x^{-3(1+2^{s})} (v^2+v) +\psi^2+\psi)
&=&-3.2^s+11
\end{eqnarray*}
car $-3.2^s+11<-3.2^s+25$ et  $-3.2^s+11<2-3.2^{s-1}$ si  $3<2^{s-1}$ \ie si $s\ge3$.
Si $s=2$, on obtient
le m\^eme r\'esultat.
En tout \'etat de cause, $v_{0,0}(x^{-3(1+2^{s})} v+\psi^2+\psi)$ est un entier impair n\'egatif.

La valuation de $x^ {-21 } (v^8+v^2)=x^ {-7 } $ en $\infty$ est
$$v_\infty(x^ {-21 } (v^8+v^2))=84-7.8=28.$$
La valuation de $(b_i x^{-3(1+2^{i})} )(v^{2^{i+1}}+v^{2^{i}}+v^2+v) $ en $\infty$ est
$$v_\infty\Big( x^{-3(1+2^{i})} (v^{2^{i+1}}+v^{2^{i}}+v^2+v) \Big)=12(1+2^{i})-7.2^{i+1}
=-2.2^{i}+12.$$
Donc la fonction $g$ a pour valuation \`a l'infini
$$v_\infty(g)=-2.2^{i}+12.$$

La valuation de $x^ {-21 } (v^8+v^2)=x^ {-7 } $ en $(0,0)$, \dots, $(0,\beta^2)$ est
$$v_{(0,0)}(x^ {-21 } (v^8+v^2))=-21+7.2=-7.$$
La valuation de $(b_i x^{-3(1+2^{i})} )(v^{2^{i+1}}+v^{2^{i}}+v^2+v) )$ en $(0,0)$ est
$$-3(1+2^{i})+7=4-3. 2^i.$$
La valuation de $g$ en $(0,0)$ est
$$v_{(0,0)}(g)=4-3. 2^s$$
si $4-3. 2^s<-7$, \ie si ${11\over3}<2^s$ \ie si $s\ge2$.

La valuation de $(v^2+v) )$ en $(0,1)$ est
$$v(v^2+v)=v\Big({v^4+v\over 1+v+v^2}\Big)=v\Big({x^7\over 1+v+v^2}\Big)=7.$$
La valuation de $(b_i x^{-3(1+2^{i})} )(v^{2^{i+1}}+v^{2^{i}}+v^2+v) )$ en $(0,1)$ est
$$v_{(0,0)}(x^{-3(1+2^{i})} )(v^{2^{i+1}}+v^{2^{i}}+v^2+v) )=-3(1+2^{i})+7=4-3. 2^i.$$

La valuation de $v^{2^{i+1}}+v^{2^{i}}+v^2+v$ en $(0,\beta)$ est
\begin{eqnarray*}
v_{(0,\beta)}(v^{2^{i+1}}+v^{2^{i}}+v^2+v)&=&
v_{(0,\beta)}((v^{2}+v)^{2^{i}}+v^2+v)\\
&=&v_{(0,\beta)}(v^2+v+1)\\
&=&7.
\end{eqnarray*}
La valuation de $(b_i x^{-3(1+2^{i})} )(v^{2^{i+1}}+v^{2^{i}}+v^2+v) )$ en $(0,\beta)$ est
$$v_{(0,\beta)}\Big((b_i x^{-3(1+2^{i})} )(v^{2^{i+1}}+v^{2^{i}}+v^2+v) \Big)
= -3(1+2^{i})+7=4-3.2^i.$$
Donc la valuation de $g$ en $(0,1)$,  $(0,\beta)$,  $(0,\beta^2)$ est
$$v_{(0,v)}(g)=4-3. 2^s$$
si $4-3. 2^s<-7$, \ie si ${11\over3}<2^s$ \ie si $s\ge2$.

    Calculons maintenant
    $$S_2=\sum_{(x,v)\in C(k)-C_\infty}\chi\pt{ g}.$$
      La valuation de $g$ en chacun de ces points finis est sup\'erieure \`a la plus faible des valuations de  
      $ (v ^2  - v ^8 )  /x ^{21}  $ et $(b_i x^{-3(1+2^{i})} )(v^{2^{i+1}}+v^{2^{i}}+v^2+v)$, elle est donc plus grande que $4-3.2^s$.
      La valuation de $g$ \`a l'infini est sup\'erieure \`a la plus faible des valuations de  
      $ (v ^2  - v ^8 )  /x ^{21}  $ et $(b_i x^{-3(1+2^{i})} )(v^{2^{i+1}}+v^{2^{i}}+v^2+v)$, elle est donc plus grande que $12-2^{s+1}$.
      Donc
      $$\deg(g)_\infty\le4(-4+3.2^s)-12+2^{s+1}=14.2^s-28.$$
      Par cons\'equent, on a
$$|S_2|\le (18-2+5+14.2^s-28) q^{1/2}=7(2^{s+1}-1)q^{1/2}.$$

Soit   $N_2$ le nombre des couples $(\alpha,v)$ tels que  $\tr\pt{{\eta} (v^2+v)}=1 $.
Alors
\begin{eqnarray*}
S_2&=&\sum_{(x,v)\in C-C_\infty} \chi\pt{ g}
=\sum_{\tr g=0} 1-N_2
=\#C-2N_2-5
\end{eqnarray*}
car $\#C_\infty=5$.
Donc 
$$\Big|N_2-{\#C\over2}\Big|={|S_2+5|\over2}\le{7\over2}(2^{s+1}-1)q^{1/2}+{5\over2}.$$

\subsubsection{Le nombre des $(\alpha,v)$ tels que $\tr\pt{{\eta} (v^2+v)}=\tr\pt{{\eta} v^3}$}

Ensuite,
nous \'evaluons le nombre des $(\alpha,v)$ tels que $\tr\pt{{\eta} (v^2+v)}=\tr\pt{{\eta} v^3}$
\ie $\tr \pt{{\eta} (v^3+v^2+v)}=0$.
Nous avons \`a calculer le nombre des $(x,v)$ tels que
$$\tr(g(x)+f(x))=0.$$
On consid\`ere la somme
$$S_3=\sum_{(x,v)\in C(k)-C_\infty}\chi\pt{ f+g}.$$
Pour v\'erifier que $f+g$ n'est pas de la forme  $\phi^2+\phi$, il suffit de calculer valuation en $(0,0)$ de $f+g+b_s\phi^2+b_s^{1/2}\phi$ comme dans la sous-section pr\'ec\'edente (\ref{2cas}). On a
$$v_{(0,0)}f=18-3.2^s\biq{et}
v_{(0,0)}(g+b_s\phi^2+b_s^{1/2}\phi)=11-3.2^s.$$
On obtient dans tous les cas une valuation impaire n\'egative. 

  Par l'analyse pr\'ec\'edente, on a
$$\deg(f+g)_\infty=21.2^s-63+14.2^s-28=35.2^s-91$$
Donc on a 
$$|S_3|\le (18-2+5+35.2^s-91)q^{1/2}=(35.2^s-70)q^{1/2}$$

Soit $N_3$  le nombre des couples $(\alpha,v)$ tels que  $\tr\pt{{\eta} (v^3+v^2+v)}=0 $.
Alors
\begin{eqnarray*}
S_3&=&\sum_{(x,v)\in C-C_\infty} \chi\pt{ f+g}
=N_3-\sum_{\tr f+g=1} 1
=2N_3-\#C+5
\end{eqnarray*}
car $\#C_\infty=5$.
 Donc 
$$\Big|N_3-{\#C\over2}\Big|={|S_3+5|\over2}\le{1\over2}(35.2^s-70)q^{1/2}+{5\over2}.$$

\subsubsection{Le nombre des $\alpha$ tels que $ X_\alpha=8q$}

Nous avons besoin d'un lemme.
\begin{lem}
\label{indep}
Soient deux fonctions $\phi$ et $\psi$ d\'efinies sur un ensemble fini $X$ \`a valeurs dans $\f2$.
Supposons que
\begin{eqnarray*}
\#\{x:\phi(x)=0\}&=&N_1 \\
\#\{x:\psi(x)=0\}&=&N_2 \\
\#\{x:\phi(x)=\psi(x)\}&=&N_3 
\end{eqnarray*}
Alors 
$$\#\{x:\phi(x)=\psi(x)=0\}= {1\over2}(N_1+N_2+N_3-N)$$
o\`u $N$ est le nombre d'\'el\'ements de $X$.
\end{lem}

{\sl D\'emonstration} --

Posons
\begin{eqnarray*}
\{x:\phi(x)=\psi(x)=0\}= N_{0,0} &\biq{,}&
\{x:\phi(x)=0,\psi(x)=1\}= N_{0,1} \\
\{x:\phi(x)=1, \psi(x)=0\}= N_{1,0} &\biq{,}&
\{x:\phi(x)=\psi(x)=1\}= N_{1,1} 
\end{eqnarray*}

On a
\begin{eqnarray*}
N_{0,0} + N_{0,1}=N_1\biq{,}
N_{0,0} + N_{1,0}=N_2\biq{,}
N_{0,0} + N_{1,1}=N_3
\end{eqnarray*}
La somme des $N_{i,j} $ \'etant \'egale \`a $N$, on a donc
$$N=\sum N_{i,j}=N_{0,0} +(N_1- N_{0,0} )+(N_2- N_{0,0}) +(N_3- N_{0,0})
=N_1+N_2+N_3- 2N_{0,0}$$
D'o\`u
$$N_{0,0}={N_1+N_2+N_3-N\over2}$$

\begin{prop}
\label{N8}
Le nombre $N$ des valeurs de $\alpha$ telles que $X_\alpha=8q$ v\'erifie
$$\Big|N-{q\over8}\Big|<23. 2^{s-1} q^{1/2}$$
pour $q\ge32$.
\end{prop}

{\sl D\'emonstration} --

D'apr\`es la proposition \ref{Xalpha}, il faut calculer  le nombre $N'$ des points $(x,v)$ tels que $\tr\pt{{\eta} v^3}=1 $ et   $\tr\pt{{\eta} (v^2+v)}=1 $.
D'apr\`es le lemme \ref{indep}, ce nombre v\'erifie
\begin{eqnarray*}
N'&=&{1\over2}(N_1+N_2+N_3-\#C) \\
&=&{1\over2}\Big(N_1-{\#C\over2}+N_2-{\#C\over2}+N_3-{\#C\over2}\Big)+{\#C\over4}
\end{eqnarray*}
et on a
\begin{eqnarray*}
\Big|N'-{\#C\over4}\Big|&=&\Big|{1\over2}\Big(N_1-{\#C\over2}+N_2-{\#C\over2}+N_3-{\#C\over2}\Big)\Big|\\
&\le&{1\over2}\Big({21.2^s-21\over2} q^{1/2}+5/2+{7\over2}(2^{s+1}-1)q^{1/2}+{5\over2}+{1\over2}(35.2^s-70)q^{1/2}+{5\over2}\Big)\\
&\le&(15/4 - 25 q^{1/2} + 91. 2^{(s-2)} q^{1/2}).
\end{eqnarray*}
Comme pour chaque $\alpha$ tel que  $\tr\pt{{\eta} v^3}=1 $ et $\tr\pt{{\eta} (v^2+v)}=1 $ il y a deux valeurs de $v$ (soit $v$ et $v+1$), le nombre $N$ de tels $\alpha$ v\'erifie donc
$$\Big|N-{\#C\over8}\Big|\le(15/8 - 25/2. q^{1/2} + 91. 2^{(s-3)} q^{1/2})$$

Comme $m$ est impair, il existe une solution de $v + v^4 = \gamma x^7$
si et seulement si la trace $\tr (\gamma x^7)$ est nulle et, dans ce cas, il y a exactement deux solutions.
Donc
	$$ \#C(k) = 2 \# \{ \tr(\gamma x^7) = 0 \} + 1 = S_7 + q + 1 $$
o\`u $S_7$ est la somme exponentielle $S_7 = \sum_{x \in k} (-1)^{\tr(\gamma x^7)}$.
Donc
	$$ | \# C(k) - q - 1 | \le 6\sqrt{q}. $$
On a
$$\Big|N-{q\over8}\Big|\le \Big|N-{\#C\over8}\Big|+\Big|{\#C\over8}-{q\over8}-{1\over8}\Big|+{1\over8}.$$
Donc le nombre $N$ v\'erifie
$$\Big|N-{q\over8}\Big|\le 15/8 - 25/2. q^{1/2} + 91. 2^{(s-3)} q^{1/2} +{3\over4}q^{1/2}+{1\over8}\le 23. 2^{s-1} q^{1/2}.$$

\section{D\'emonstration de l'\'evaluation de $\|\ch f\|_4^4$ (proposition \ref{eval})}
\label{demeval}

On d\'eduit facilement de la proposition \ref{N8} le calcul de la valeur de $\|\ch f\|_4^4$.
Sachant que
\begin{eqnarray*}
\Big|N-{q\over8}\Big|&\le&23.2^{s-1} q^{1/2}\\
\Big|N_0-{q\over2}\Big|&\le&3q^{1/2}+1
\end{eqnarray*}
calculons
\begin{eqnarray*}
\|\ch f\|_4^4 &=&q^2+\sum_{\textstyle{\alpha\ne0\atop \alpha\in V_m}}X_\alpha\\
&=&3q^2+8q(N-q/8)+2q(N_0-q/2).
\end{eqnarray*}
D'o\`u
\begin{eqnarray*}
|\|\ch f\|_4^4 -3q^2|
&\le& 8q\Big|N-{q\over8}\Big|+2q\Big|N_0-{q\over2}\Big|\\
&\le& 185.2^{s-1}q^{3/2}.
\end{eqnarray*}

\section{D\'emonstration des bornes de $\|\ch f\|_\infty$ (propositions \ref{borneinf} et \ref{bornesup})}
\label{dembornes}

\subsection{Borne inf\'erieure}

L'\'evaluation du nombre des $\alpha$ tels que $\tr \ell=1$ dans la proposition \ref{Xalpha} donne:
$$2q^2-6q^{3/2}\le \|\ch f\|_4^4.$$

On a
$$\sum_{\alpha\in k^*} X_\alpha\ge 2q N\ge 2q {q-6\sqrt q\over2}=q^2-6q^{3/2}$$
et
$$\|\ch f\|_4^4=q^2+\sum_{\alpha\in k^*} X_\alpha\ge 2q^2-6q^{3/2}.$$

Comme il est facile de montrer que
$$\|\ch f\|_4^4\le q\|\ch f\|_\infty^2$$
nous obtenons
$2q-6q^{1/2}\le \|\ch f\|_\infty^2$, donc $\sqrt{2q}-3\sqrt 2\le \|\ch f\|_\infty$,
d'o\`u le r\'esultat si $m\ge7$, parce que $\|\ch f\|_\infty$ est un entier divisible par $2^{\lceil m/3\rceil}$.
Le r\'esultat pour $m\le7$ est connu (cf. remarque \ref{connu}).

On a, de plus
$$\|\ch f\|_4^4\ge 3q^2-185.2^{s-1}q^{3/2}$$
par la proposition \ref{eval}.
On en d\'eduit que
pour que $\|\ch f\|_4^4$ d\'epasse $2q^2$, il suffit que $m\ge15+2s$.
Pour des raisons de divisibilit\'e, $\|\ch f\|_\infty$ est alors plus grand que $\sqrt{2q}+ 2^{\lceil{m\over 3}\rceil}$.

\subsection{Borne sup\'erieure}

On a, d'apr\`es la borne de Weil
$$|\ch f(v)|= \Bigl| \sum_{x\in V_m}\chi{(f(x)+v\cdot x)}\Bigr|\le6\sqrt q.$$


\begin{thebibliography}{XXXX}

\bibitem{ca}{Anne Canteaut }
{\it Differential cryptanalysis of Feistel ciphers
and differentially  $\delta$-uniform mappings,}
in Selected Areas on Cryptography, SAC'97, pages 172-184, Ottawa, Canada, 1997.

\bibitem{cccf}{A.~Canteaut,  C.~Carlet, 
P.~Charpin,  C.~Fontaine}
{\it Propagation characteristics et correlation-immunity of
highly nonlinear Boolean functions,}
Advances in cryptology, EUROCRYPT 2000 (Bruges), 507--522,
Lecture Notes in Comput. Sci., Vol. 1807, Springer, Berlin, 2000.

\bibitem{brv}{Barth, Rolland, V\'eron}
{\it Cryptographie}, Herm\`es, Paris, 2005.

\bibitem{bo}{E. Bombieri,} 
{\it On exponential sums in finite fields.}
Amer. J. Math., 88, 1966, pp. 71-105.

\bibitem{ca2}{C. Carlet},
{\it On cryptographic complexity of Boolean
functions},
Proceedings of the Sixth Conference on Finite Fields with 
Applications to Coding Theory, Cryptography et Related Areas
(G.L. Mullen, H. Stichtenoth et H. Tapia-Recillas
Eds),
Springer (2002) pp. 53-69.


\bibitem{ca3}{C. Carlet},
{\it On the algebraic thickness et 
non-normality  of Boolean functions, with developments on symmetric 
functions}, submitted to IEEE Trans. Inform. Theory.

\bibitem{cal}
{G.~Cohen, 
I.~Honkala, S.~Litsyn, 
A.~Lobstein},
{\it Covering codes}.
North-Holland Mathematical Library, 54,
North-Holland Publishing Co., Amsterdam (1997). 


\bibitem{cf}
{C. Fontaine},
{\it Contribution \`a la recherche de fonctions bool\'eennes
hautement non  lin\'eaires
     et au marquage d'images en vue de la protection des 
droits d'auteur},
Th\`ese, Universit\'e Paris VI (1998).

\bibitem{ho}
{X. Hou,}
{Covering radius of the Reed-Muller code $R(1,7)$---a simpler  proof.}  J. Combin. Theory Ser. A  74  (1996),  no. 2, 337--341.
                
 \bibitem{le}
 {Leander, Nils Gregor} 
 {Monomial bent functions.}  
 IEEE Trans. Inform. Theory  52  (2006),  no. 2, 738--743.

 \bibitem{lz}
{Langevin, P.; Zanotti, J.-P.}
 {Nonlinearity of some invariant Boolean functions.}  
 Des. Codes Cryptogr.  36  (2005),  no. 2, 131--146. 

\bibitem{ma}
{Sel\c{c}uk Kavut, Subhamoy Maitra and Melek D. Y{\"u}cel}
{There exist Boolean functions on $n$ (odd) variables having nonlinearity $> 2^{n-1} - 2^{\frac{n-1}{2}}$ if and only if $n > 7$}, pr\'epublication,
{\tt http://eprint.iacr.org/2006/181}


\bibitem{mm}
{C. Moreno et O. Moreno}
{The MacWilliams-Sloane conjecture on the tightness of the  Carlitz-Uchiyama bound and the weights of duals of BCH codes.}  IEEE Trans. Inform. Theory  40  (1994),  no. 6, 1894--1907.

\bibitem{mn}
{Daniel Maisner et Enric Nart},
{\it Zeta functions of supersingular curves
of genus 2},
arXiv:math.NT/0408383



\bibitem{ms}
{F.J. MacWilliams et N.J.A. Sloane},
{\it The Theory
of Error-Correcting Codes,} North-Holland, Amsterdam (1977).



\bibitem{pw} {N. Patterson et D. Wiedemann},
{\it  The covering
radius of the $(2\sp{15},\,16)$ Reed-Muller code is at least 
$16\,276$},
IEEE Trans. Inform. Theory 29, no. 3 (1983), 354-356.


\bibitem{ro1}
{F. Rodier},
{\it Sur la non-lin\'earit\'e des fonctions bool\'eennes},
 Acta
Arithmetica, vol 115,  (2004), 1-22, 
preprint:
arXiv: math.NT/0306395.


\bibitem{ro2}
{F. Rodier},
{\it On the nonlinearity of Boolean functions}, Proceedings of WCC2003,
Workshop on coding et cryptography 2003   (D.~Augot,
P.~Charpin, G.~Kabatianski eds), INRIA (2003), pp. 397-405. 

\bibitem{se}
{J-P. Serre,} 
{\it Majorations de sommes exponentielles.}
Journ\'ees Arithm\'etiques de Caen (Univ. Caen, Caen, 1976), pp. 111-126. Ast\'erisque No. 41-42,
Soc. Math.France, Paris, 1977.

\bibitem{sti}
{ H. Stichtenoth,} {Algebraic Function Fields et Codes,}
Springer, 1993.


\bibitem{st}
{P. St\u anic\u a}, 
{\it Nonlinearity, local et global avalanche characteristics of
balanced Boolean functions}, 
Discrete Math. 248 (2002), no. 1-3, 181--193.

\bibitem{gv1}
 {G. van der Geer, M. van der Vlugt}, {\it Reed-Muller codes and supersingular curves. I}, Compositio Math. 84, (1992), 333-367.


\bibitem{gv2}
{G. van der Geer, M. van der Vlugt}, {\it Supersingular Curves of
Genus 2 over finite fields of Characteristic 2} , Math. Nachr. 159, (1992),
73-81.

\bibitem{zz}{Xian-Mo Zhang and
Yuliang Zheng},
{\it GAC |the Criterion for Global Avalanche
Characteristics of Cryptographic Functions},
Journal of Universal Computer Science, vol. 1, no. 5 (1995), 316-333


\end{thebibliography}
\end{document}